\documentclass{article}

\makeatletter
\renewcommand{\@biblabel}[1]{[#1]\hfill}
\makeatother

\let\OLDthebibliography\thebibliography
\renewcommand\thebibliography[1]{
  \OLDthebibliography{#1}
  \setlength{\parskip}{0pt}
  \setlength{\itemsep}{3pt plus 0.3ex}
}

\makeatletter
\def\blankfootnote{\xdef\@thefnmark{}\@footnotetext}
\makeatother

\usepackage[T1]{fontenc}

\usepackage[letterpaper, portrait, left=1in, right=1in, top=0.9in, bottom=1in, footnotesep=1.5\baselineskip]{geometry}
\usepackage[dvipsnames]{xcolor}
\usepackage{contour}
    \contourlength{0.5pt}
\usepackage{tocloft}

\usepackage{environ} 
\usepackage[framemethod=tikz]{mdframed}
\usepackage{wrapfig}
\usepackage{mathtools}
    \usetikzlibrary{calc, math, cd, arrows.meta, braids, decorations.pathreplacing, decorations.markings, decorations.pathmorphing, bending, shapes}
    \usepackage{tikz-3dplot}
    \tikzset{vertex/.style={draw, shape=circle, inner sep=1.5pt, minimum size=4pt}}
    \tikzset{<->/.tip={Latex}}
    \tikzset{shorten > = 2pt, shorten <=2pt}
    \tikzset{smallnode/.style={every node/.style={draw, fill=black, shape=circle, inner sep=0pt, minimum size=4pt}, scale=0.7}}
    \tikzset{drawnode/.style={fill,shape=circle,inner sep=0pt, minimum size=3pt}}
\usepackage{pgfornament}
\usepackage{fourier-orns} 
\usepackage{stmaryrd}
\usepackage{halloweenmath}
\usepackage[shortlabels]{enumitem}
    \newlength{\circlabelwidth}
    \setlist{nosep}
    \setlist[enumerate]{label=\textup{\arabic*.}}
    \newlist{subprob}{enumerate}{2}
        \setlist[subprob,1]{label={(\roman*)}}
        \setlist[subprob,2]{label={(\arabic*)}}
    \setlist[itemize]{labelindent=10pt,labelwidth=\circlabelwidth,leftmargin=!,label=$\circ$}
    \newlist{problems}{enumerate}{3}
        \setlist[problems,1]{before=\setupstar,label=\textup{\arabic*.}, itemsep=2pt, topsep=8pt,ref=\textup{\arabic*}}
        \setlist[problems,2]{before=\setupstar,label=(\alph*),parsep=0pt}
        \setlist[problems,3]{before=\setupstar,label=(\roman*),parsep=0pt}

\usepackage[normalproof]{dopestyle}

\makeatletter
    \renewcommand\@makefntext[1]{\leftskip=0em\hskip-0em\@makefnmark\,#1}
\makeatother

\surroundwithmdframed[skipabove=0.5\baselineskip, skipbelow=0.5\baselineskip, leftmargin=3pt, rightmargin=3pt, innerleftmargin=7pt, innerrightmargin=7pt, roundcorner=10pt, linewidth=2pt, linecolor=red!40, backgroundcolor=red!5]{headsup}

\surroundwithmdframed[topline=false, bottomline=false, innertopmargin=0pt, innerbottommargin=0pt, innerleftmargin=2pt, innerrightmargin=2pt, linewidth=0.2mm]{quote}

\surroundwithmdframed[topline=false, bottomline=false, innertopmargin=2.5pt, innerbottommargin=2.5pt, innerleftmargin=-10pt, leftmargin=-10pt, innerrightmargin=-10pt, rightmargin=7.7pt, linewidth=0.4mm]{quote}

\NewEnviron{method}{
    \parbox{\textwidth}{
        \textbf{\textsc{Method}}
        \begin{mdframed}[innerleftmargin=4pt,innerrightmargin=4pt,skipabove=3pt,skipbelow=0pt]
            \BODY
        \end{mdframed}
    }
}

\theoremstyle{itcaps}
\newtheorem{theorem}{Theorem}
\newtheorem{corollary}[theorem]{Corollary}

\theoremstyle{alert}

\theoremstyle{solved}

\theoremstyle{caps}

\newtheorem{exampleprimitive}[theorem]{Example}

    \crefname{exercise}{Exercise}{Exercises}

\theoremstyle{remark}

\numberwithin{equation}{section}

\newcommand*{\newword}[2][]{\emph{#2}\index{%
    \ifx&#1&%
       #2%
    \else%
       #1%
    \fi}%
} 
\newcommand*{\oldword}[2][]{#2\index{%
    \ifx&#1&%
       #2%
    \else%
       #1%
    \fi}%
} 

\DeclareMathOperator{\conv}{conv}

\let\phi\varphi

\let\epsilon\varepsilon

\let\oldchi\chi
\renewcommand{\chi}{\raisebox{1pt}{$\oldchi$}}

\newcommand{\inserttitle}{Layer number of the grid}
\title{\inserttitle}
\date{}

\allowdisplaybreaks

\usepackage{caption}
\begin{document}

\setlength{\abovedisplayskip}{6pt plus 2pt minus 4pt}
\setlength{\abovedisplayshortskip}{1pt plus 3pt}
\setlength{\belowdisplayskip}{6pt plus 2pt minus 4pt}
\setlength{\belowdisplayshortskip}{6pt plus 2pt minus 2pt}

\setlength{\parindent}{10pt}

{\centering
{\LARGE \scshape explicit bounds for the\\ layer number of the grid\par\vspace{0.5\baselineskip}}
{\scalebox{1.1}{Travis Dillon and Narmada Varadarajan}\par}
}\vspace{\baselineskip}

\begin{abstract}
    The number of steps required to exhaust a point set by iteratively removing the vertices of its convex hull is called the \emph{layer number} of the point set. This article presents a short proof that the layer number of the grid $\{1,2,\dots,n\}^d$ is at most $\frac{1}{4}dn^2+1$, significantly improving the dependence on $d$ in the best-known upper bound. We also prove a lower bound of $\frac{1}{2}d(n-1)+1$, which shows that the layer number of the grid is linear in $d$.
    \blankfootnote{\null\\[-2.2em]\\
        \noindent Dillon was supported by a National Science Foundation Graduate Research Fellowship under Grant No.\ 2141064.
    }
\end{abstract}

\section{Introduction}

The convex layers of a point set in $\R^n$ are constructed by selecting the vertices of its convex hull, removing them from the set, and repeating. More formally, suppose that $X$ is a finite point set in $\R^n$. Setting $C_1 = X$, we inductively define $L_i$ as the vertices of the convex hull of $C_i$ and set $C_{i+1} = C_i\setminus L_i$.
The greatest integer $k$ for which $C_k \neq \emptyset$ is called the \emph{layer number} of $X$ and is denoted $L(X)$. This process is sometimes called \emph{peeling}, and the set of layers is whimsically referred to as the \emph{onion}.

Convex layers were introduced by Barnett as one of several possible methods of ordering high-dimensional data sets \cite{first-ordering}. The peeling process was studied algorithmically by Eddy \cite{eddy} and Chazelle \cite{chazelle}, and has since found applications in outlier detection \cite{outlier-detection}, recognition of projectively deformed point sets \cite{layer-deform-recognition}, and fingerprint recognition \cite{fingerprint}.

More recently, mathematical attention has turned toward calculating the layer number of various classes of point sets. In 2004, Dalal proved that the expected number of layers in an random $n$-point set distributed uniformly in a $d$-dimensional ball is $\Theta_d(n^{2/d+1})$ \cite{dalal}. (We use subscripts in order-of-magnitude notation to denote the variables on which the hidden constants depend.) Following this work, Choi, Joo, and Kim studied the layer number of point sets that avoid clustering \cite{alpha}, and they prove that $L(X) = O_{\alpha,d}\big(|X|^{(d+1)/2d}\big)$ for any $\alpha$-evenly distributed point set in the unit ball of $\R^d$. (The parameter $\alpha$ is a measure of the evenness of spread.) Ambrus, Nielsen, and Wilson improved this result, showing that any such set has $O_{\alpha,d}\big(|X|^{2/d})$ layers when $d \geq 3$ \cite{alpha2}.

Given two integers $a \leq b$, we write $[a,b] = \{a,a+1,\dots,b-1,b\}$ for the integer range from $a$ to $b$, and we set $[n] = [1,n]$. The problem of determining the layer number of the integer grid $[n]^d$ has attracted attention as a notable special case of the layer number problem \cite{alpha2,alpha}. In 2013, Har-Peled and Lidick\'y proved that $L([n]^2) = \Theta(n^{4/3})$ as $n\to\infty$, solving the two-dimensional problem \cite{two-dim-grid}. Ambrus, Hsu, Peng, and Jan \cite{higher-grid} proved that $L([n]^d) = \Omega_d\big(n^{2d/(d+1)}\big)$ by extending Har-Peled and Lidik\'y's proof, and they conjecture that this is tight. Since the grid $[n]^d$ is evenly distributed, Ambrus--Nielsen--Wilson's upper bound \cite{alpha2} shows that $L([n]^d) = O_d(n^2)$, where the hidden constant is exponential in $d$.

In this paper, we present a simple argument that significantly improves the dependence on $d$:

\begin{theorem}\label{thm:upper}
    $L([2n+1]^d) \leq dn^2 + 1$.
\end{theorem}

\cref{thm:upper} extends to a bound on $L([n]^d)$ for every $n$ using the fact that $L(X)\leq L(Y)$ whenever $X \subseteq Y\hspace{-0.15em}$. This can be proven by showing that the containment continues to hold at each step of the peeling process (see \cite[Lemma 2.1]{alpha}, for example). We also prove a new lower bound on the layer number of the grid:

\begin{theorem}\label{thm:lower}
    $L([2n+1]^d) \geq dn+1$.
\end{theorem}

This lower bound is weaker than current lower bounds in the setting where the dimension is fixed and $n\to\infty$. However, in the opposing regime where $n$ is fixed and $d\to\infty$, \cref{thm:upper,thm:lower} completely determine the order of magnitude of the layer number.

\begin{corollary}
    $L([n]^d) = \Theta_n(d)$.
\end{corollary}

\section{Proofs}

In the following proofs, we consider the translated point set $[-n,n]^d$ in place of $[2n+1]^d$. The advantage of this is that the group of symmetries $\Gamma$ of $[-n,n]^d$ is generated by permutations and negations of coordinates. Moreover, the layers of $[-n,n]^d$ are invariant under the action of $\Gamma$: If a point $y$ is peeled at layer $i$, so are all of its images under $\Gamma$. This fact is especially important for the proof of \cref{thm:lower}.

\begin{proof}[Proof of \cref{thm:upper}]
    Let $L_i$ denote the $i$th layer of the peeling process and $C_i = \conv(L_i) \cap [-n,n]^d$ denote the set of grid points remaining after $i-1$ peeling steps. For each $i$, pick an integer point $y_i \in L_i$ of maximum norm and set $r_i = \|y_i\|$. Since each point in $[-n,n]^d$ has integer coordinates, $r_i^2$ is an integer. Any point $y \in L_i$ of norm $r_i$ is an extreme point of $C_i$, since the tangent hyperplane at $y$ to the sphere of radius $r_i$ contains all points of $C_i$ in one of its half-spaces. This means that every point of $L_{i+1}$ has norm strictly less than $r_i$, so $r_{i+1} < r_i$.
    
    Therefore $dn^2 = r_1^2 > r_2^2 > \cdots > r_{L([-n,n]^d)}^2$ is a strictly decreasing nonnegative set of integers, so it has at most $dn^2 + 1$ entries.
\end{proof}

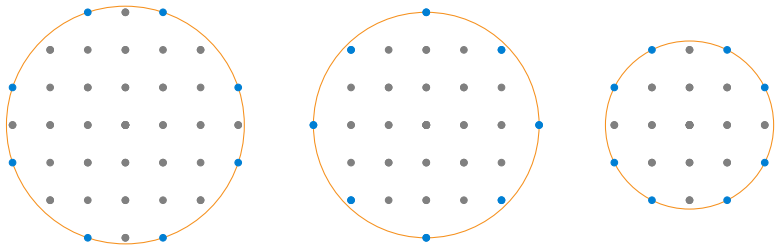
\begin{figure}[b]
\centering
\definecolor{lightblue}{HTML}{007FD4}
\newcommand{\vertexcolor}{lightblue}
\newcommand{\circlecolor}{BurntOrange}
{\centering
\begin{tikzpicture}[every node/.style={draw=gray,shape=circle,fill=gray,inner sep=0pt, minimum size=2.5pt}, scale=0.5]
    \begin{scope}[shift={(9,0)}]
        \draw[\circlecolor] (0,0) circle[radius=sqrt(10)];
        \foreach \x/\y in {0/0, 0/1, 0/2, 1/1, 1/2, 2/2, 0/3}{
            \foreach \s in {+,-}{
                \foreach \t in {+,-}{
                    \node at (\s\x,\t\y) {};
                    \node at (\s\y,\t\x) {};
                }
            }
        }
        \foreach \x/\y in {1/3}{
            \foreach \s in {+,-}{
                \foreach \t in {+,-}{
                    \node[\vertexcolor] at (\s\x,\t\y) {};
                    \node[\vertexcolor] at (\s\y,\t\x) {};
                }
            }
        }
    \end{scope}
    \begin{scope}[shift={(17,0)}]
        \draw[\circlecolor] (0,0) circle[radius=sqrt(9)];
        \foreach \x/\y in {0/0, 0/1, 0/2, 1/1, 1/2}{
            \foreach \s in {+,-}{
                \foreach \t in {+,-}{
                    \node at (\s\x,\t\y) {};
                    \node at (\s\y,\t\x) {};
                }
            }
        }
        \foreach \x/\y in {0/3, 2/2}{
            \foreach \s in {+,-}{
                \foreach \t in {+,-}{
                    \node[\vertexcolor] at (\s\x,\t\y) {};
                    \node[\vertexcolor] at (\s\y,\t\x) {};
                }
            }
        }
    \end{scope}
    \begin{scope}[shift={(24,0)}]
        \draw[\circlecolor] (0,0) circle[radius=sqrt(5)];
        \foreach \x/\y in {0/0, 0/1, 0/2, 1/1}{
            \foreach \s in {+,-}{
                \foreach \t in {+,-}{
                    \node at (\s\x,\t\y) {};
                    \node at (\s\y,\t\x) {};
                }
            }
        }
        \foreach \x/\y in {1/2}{
            \foreach \s in {+,-}{
                \foreach \t in {+,-}{
                    \node[\vertexcolor] at (\s\x,\t\y) {};
                    \node[\vertexcolor] at (\s\y,\t\x) {};
                }
            }
        }
    \end{scope}
\end{tikzpicture}\par}
\caption{Steps 3, 4, 5 in the peeling process of $[7]^2$. The blue points are the vertices of the convex hull at each step, and the orange circle is the sphere of radius $r_i$ in the proof of \cref{thm:upper}.}
\label{fig:circles}
\end{figure}

\cref{fig:circles} provides a visual interpretation of the argument in the proof of \cref{thm:upper}.

\begin{proof}[Proof of \cref{thm:lower}]
    The key tool is a relation $\prec$ on the points of the grid so that $x \prec y$ implies that $x$ is peeled after $y$. The length of any chain with respect to this relation is a lower bound on the number of layers.

    Given $x,y \in [-n,n]^d$, we write $x \prec y$ if there is a value $k$ so that $x_i = y_i$ for every $i \neq k$ and $|x_k| < |y_k|$. We now show that $x$ is peeled after $y$ whenever $x \prec y$. Using the symmetry group, we may assume that all coordinates of $x$ and $y$ are nonnegative and that $k=1$. The symmetric image $y' = (-y_1, y_2,\dots,y_d)$ of $y$ is in the same layer as $y$ and 
    \[
        x = \frac{1}{2} \Big(1+\frac{x_1}{y_1}\Big) y + \frac{1}{2} \Big(1-\frac{x_1}{y_1}\Big) y'.
    \]
    Since $x$ is a nontrivial convex combination of two vertices, it is not an extreme point in the same layer as $y$ and must be peeled later.

    Finally, form a chain $(0,0,\dots,0) = x^1 \prec x^2 \prec \cdots \prec x^{nd+1} = (n,n,\dots,n)$ by increasing exactly one coordinate by 1 at each step. We conclude that $L([-n,n]^d) \geq nd+1$.
\end{proof}\vspace{1\baselineskip}

\noindent
{\Large\scshape acknowledgements}\\[1ex]
We thank Shubhangi Saraf and Gergely Ambrus for helpful conversations during the preparation of this paper. We are also very grateful to Callie Garst for originally connecting the two authors.
\vspace{2\baselineskip}

\renewcommand{\refname}{references}
\bibliographystyle{amsplain-nodash}
\bibliography{bibliography}

\ifx\undefined\bysame
\newcommand{\bysame}{\leavevmode\hbox to3em{\hrulefill}\,}
\fi
\begin{thebibliography}{10}

\bibitem{higher-grid}
Gergely Ambrus, Alexander Hsu, Bo~Peng, and Shiyu Jan, {\em The layer number of
  grids}, 2020, \href{https://arxiv.org/abs/2009.13130}{arXiv:2009.13130}
  [math.MG].

\bibitem{alpha2}
Gergely Ambrus, Peter Nielsen, and Caledonia Wilson, {\em New estimates for
  convex layer numbers}, Discrete Mathematics {\bf 344} (2021), 112424.

\bibitem{first-ordering}
V.~Barnett, {\em The ordering of multivariate data}, Journal of the Royal
  Statistical Society {\bf 139} (1976), 318--355.

\bibitem{chazelle}
Bernard Chazelle, {\em On the convex layers of a planar set}, IEEE Transactions
  on Information Theory {\bf 31} (1985), 509--517.

\bibitem{alpha}
Ilkyoo Choi, Weonyoung Joo, and Minki Kim, {\em The layer number of
  $\alpha$-evenly distributed point sets}, Discrete Mathematics {\bf 343}
  (2020), 112029.

\bibitem{dalal}
Ketan Dalal, {\em Counting the onion}, Random Structures and Algorithms {\bf
  24} (2004), 155--165.

\bibitem{eddy}
W.~F. Eddy, {\em Convex hull peeling}, COMPSTAT 5th Symposium, 1982, 42--47.

\bibitem{two-dim-grid}
Sariel Har-Peled and Bernard Lidick\'y, {\em Peeling the grid}, SIAM Journal on
  Discrete Mathematics {\bf 27} (2013), 650--655.

\bibitem{outlier-detection}
Victoria~J. Hodge and Jim Austin, {\em A survey of outlier detection
  methodologies}, Artificial Intelligence Review {\bf 22} (2004), 85--126.

\bibitem{fingerprint}
Marios Poulos, Sozon Papavlasopoulos, and Vasilios Chrissikopoulos, {\em An
  application of the onion peeling algorithm for fingerprint verification
  purposes}, Journal of Information and Optimization Sciences {\bf 26} (2005),
  665--681.

\bibitem{layer-deform-recognition}
Tom\'a\v{s} Suk and Jan Flusser, {\em Convex layers: A new tool for recognition
  of projectively deformed point sets}, Computer Analysis of Images and
  Patterns, 1999, 454--461.

\end{thebibliography}

\vspace{1.2\baselineskip}

\noindent
\textsc{Department of Mathematics, Massachusetts Institute of Technology, Cambridge, MA, USA}\\
\textit{email:} \texttt{dillont@mit.edu}\vspace{\baselineskip}

\noindent\textsc{Department of Mathematics, University of Toronto, Toronto, ON, Canada}\\
\textit{email:} \texttt{narmada.varadarajan@mail.utoronto.ca}
\end{document}